# GENERIC-CASE COMPLEXITY, DECISION PROBLEMS IN GROUP THEORY AND RANDOM WALKS

ILYA KAPOVICH, ALEXEI MYASNIKOV, PAUL SCHUPP, AND VLADIMIR SHPILRAIN

ABSTRACT. We give a precise definition of "generic-case complexity" and show that for a very large class of finitely generated groups the classical decision problems of group theory - the word, conjugacy and membership problems - all have linear-time generic-case complexity. We prove such theorems by using the theory of random walks on regular graphs.

## Contents

## 1. MOTIVATION

Algorithmic problems such as the word, conjugacy and membership problems have played an important role in group theory since the work of M. Dehn in the early 1900's. These problems are "decision problems" which ask for a "yes-or-no" answer to a specific question. For example, the word problem for a finitely presented group $G = \langle x_1, \ldots, x_k \, | \, r_1, \ldots, r_m \rangle$ asks, given a word $w$ in $\{x_1, \ldots, х_k\}^{\pm 1}$, whether or not this word represents the identity element of $G$. The classical result of P. S. Novikov and of W. Boone states that there exists a finitely presented group with unsolvable word problem. This implies that most other problems (the conjugacy, membership, isomorphism, and order problems) are also unsolvable in the class of all finitely presented groups (see the survey papers [1, 48] for a detailed exposition).

With the advance of modern computers interest in algorithmic mathematics shifted to the realm of decidable problems with a particular emphasis on complexity of algorithms, and in the 1970s modern complexity theory was born. It quickly turned out that some decidable problems which one would really like to solve are too difficult to be solved in full generality on actual computers. Among different possible complexity measures the most important for us here is *time complexity*. Usually, algorithms with linear, or quadratic, or sometimes even with high degree polynomial time complexity, are viewed as fast algorithms. Fortunately, several classes of infinite groups have fast algorithms for their decision problems. For example, the word and conjugacy problems for any word-hyperbolic group are solvable in linear and in quadratic time respectively, and the word problem for a linear group over the field of rational numbers can be solved in cubic time. On the other hand, there are finitely presented groups whose word problem has arbitrarily high time complexity. For a group with exponential time complexity of the word problem any algorithm solving the word problem needs at least exponentially many steps (in length of the word) to halt on infinitely many inputs. This type of analysis concerns the worst-case behavior of an algorithm and is now often called *worst-case complexity*.

Many algorithms for solving the word problem in finitely presented groups are difficult to analyze and their worst-case complexity is not known. For example, for the Magnus algorithm for the word problem for one-relator groups [45] we do not even know if the complexity is bounded above by any fixed tower of exponentials. Yet anyone who has conducted computer experiments with finitely presented groups knows that there is often some kind of an easy "fast check" algorithm which quickly produces a solution for "most" inputs of the problem. This is true even if the worst-case complexity of the particular problem is very high or the problem is unsolvable. Thus many group-theoretic decision problems have a very large set of inputs where the (usually negative) answer can be obtained easily and quickly. Indeed, our intuition on the subject has been formed by computer experiments and the main purpose of this paper is to explain some of this phenomenon. It turns out that a precise mathematical explanation comes from the theory of random walks on regular graphs.

The kind of situation which we have in mind is often analogous to the use of Dantzig's Simplex Algorithm for linear programming problems. This algorithm is used hundreds of times daily and in practice almost always works quickly. The examples of V. Klee and G. Minty [41] showing that one can make the simplex algorithm take exponential time are very special. A "generic" or "random" linear programming problem is not "special", and the algorithm works quickly. Observations of this type led to the development of *average-case complexity*. There are several different approaches to the average-case complexity, but they all involve computing the expected value of the running time of an algorithm with respect to some measure on the set of inputs (for example, see [37, 42]).

To study *generic-case* complexity, which deals with the performance of an algorithm on "most" inputs, we first need a notion of which sets are *generic*. Let $\nu$ be a probability distribution on $X^*$, or, more generally, an arbitrary additive function with values in $[0, 1]$ defined on some subsets of the set $X^*$ of all finite words over a finite alphabet $X$. A subset $T \subset X^*$ is called *generic with respect to* $\nu$ if $\nu(X^* - T) = 0$. Then, for example, we would say that an algorithm $\Omega$ has *polynomial-time generic-case complexity with respect to* $\nu$, if $\Omega$ runs in polynomial time on all inputs from some subset $T$ of $X^*$ which is generic with respect to $\nu$. Of course, we can define generic-case complexity being in any complexity class $\mathcal{C}$, not only for polynomial time.

Thus "generic-case" complexity is in the spirit of but quite different from average-case complexity in several respects. First of all, in average-case complexity the decision problem considered must be decidable and one has to have a total algorithm to solve it. One is then interested in the expected value of the running time of the algorithm. On the other hand, in generic-case complexity we consider the behavior of the algorithm only on a generic set $T$ and completely ignore its behavior elsewhere. Thus we consider partial algorithms which may only halt on the set $T$ and the total problem being considered can have arbitrarily high worst-case complexity or even be undecidable.

The general idea of generic behavior in the context of group theory was introduced by M. Gromov [34] when he defined the class of word-hyperbolic groups. Gromov indicated that "most" finitely presented groups are word-hyperbolic. This was made precise by A. Ol'shanskii [50] and also by C. Champetier [19] who formalized the notion of a "generic" group-theoretic property. Further research on generic group-theoretic properties has been done by C. Champetier [19, 20, 21], G. Arzhantseva [7, 8, 9, 6], A. Zuk (unpublished), P.-A. Cherix with co-authors [22, 23] and others. Recently M. Gromov [36] pushed his ideas about "random groups" further with the goal of constructing finitely presentable groups that do not admit uniform embeddings into a Hilbert space.

The notion of genericity in the work cited above concerns the collection of all finitely presented groups. In this paper we shift the focus to considering generic properties of algorithmic problems in *individual* groups with respect to *asymptotic density* (see Section 3).

We are grateful to Laurent Bartholdi and Tatiana Smirnova-Nagnibeda for many helpful discussions regarding random walks on groups and graphs, to Carl Jockusch and Frank Stephan for discussions on the general idea of generic-case complexity and to Bogdan Peternko for suggesting the use of Stolz' theorem.

2. Algorithms and decision problems for groups

**Convention 2.1.** We follow the book *Computational Complexity* of C. Papadimtriou [51] for our conventions on computational complexity. Recall that a *complexity class* is determined by specifying a *model of computation* (which for us is always a multi-tape Turing machine), a *mode of computation* (e.g. deterministic or non-deterministic), *resources* to be controlled (e.g. time and space) and *bounds* for each controlled resource, that is functions $f(x)$ such that for each input word $w$ at most $f(|w|)$ units of the specified resource needs to be extended by an appropriate Turing machine to reach a decision.

In this paper, unless specified otherwise, when talking about a "complexity class $\mathcal{C}$", we assume that the resources to be controlled in the definition of $\mathcal{C}$ are either time or space or a combination thereof. We also assume that the collection of functions bounding each resource consists of *proper complexity functions* $f(n) > 0$ (see [51] and Section Section:lang below) and that for any function $f(n)$ in this collection and for any integer $C \geq 1$ the function $Cf(Cn+C)+C$ also belongs to this collection. Most reasonable complexity classes, such as linear time, polynomial space, log-space, etc, satisfy these restrictions.

Recall that a *decision problem* is a subset $\mathcal{D}$ of the set
$$(X^*)^k = X^* \times \cdots \times X^*$$
($k \geq 1$ factors), where $X^*$ is the set of all words on a finite alphabet $X$. (By introducing an extra alphabet symbol "," we could view a $k$-tuple of words $(w_1, w_2, \ldots, w_k) \in (X^*)^k$ as a single word in the alphabet $X \cup \{,\}$.)

In this section we focus on three classical decision problems for a given finitely generated group $G$: the *word problem* (WP), the *conjugacy problem* (CP), and the *subgroup membership problem* (MP). (Our approach is quite general and can be applied to other group-theoretic decision problems, such as the order of an element problem.) To formulate these problems precisely one needs to specify exactly how the group $G$ is "given". To do this, one chooses a finite set of generators $A$ of a group $G$, that is, one fixes a map $\pi : A \to G$ such that $G = \langle \pi(A) \rangle$. To simplify notation we identify elements of $A$ with their images under $\pi$ in $G$. Put $X = A \cup A^{-1}$. Thus every word $w \in X^*$ represents an element $\pi(w) \in G$.

Now we are ready to formulate the algorithmic problems above *with respect to the given set of generators $A$*:

(WP) Given a word $w \in X^*$ determine whether or not $w$ represents the identity element in $G$ (symbolically, $w =_G 1$). Thus
$$WP(G, A) := \{w \in X^* \,|\, w =_G 1\}.$$

(CP) Given two words $u, v \in X^*$ determine whether they represent conjugate elements of $G$ or not. Thus
$$CP(G, A) := \{(u, v) \in X^* \times X^* \,|\, \pi(u), \pi(v) \text{ are conjugate in } G\}.$$

(MP) Let $H \leq G$ be a fixed finitely generated subgroup. Given a word $u \in X^*$ determine whether or not $u$ belongs to $H$. Thus
$$MP(G, H, A) := \{w \in X^* \,|\, \pi(w) \in H\}.$$

**Convention 2.2.** We call these problems the *A-versions* of the corresponding problem about $G$ to emphasize the choice of generators $A$. We use the notation **D** to denote a problem about a group $G$ and we denote by $\mathcal{D}_A$ the $A$-version of **D** corresponding to the finite generating set $A$ of $G$. Thus if **D** is the word problem for $G$, then $\mathcal{D}_A = WP(G, A)$. If **D** is a problem about a group $G$ and $\mathcal{C}$ is a complexity class, we say that **D** is *solvable* for $G$ with complexity in $\mathcal{C}$ if for every finite generating set $A$ of $G$ the language $\mathcal{D}_A$ is in $\mathcal{C}$.

Of course, instead of the problems over $X^*$ one can consider decision problems only over freely reduced words, that is, decision problems $\mathcal{D} \subset F(A)^k$, where $F(A)$ is the free group on $A$. Since one can easily (in linear time) reduce a word in $X^*$ to its reduced form in $F(A)$ these two decision problems are equivalent with respect to time complexity classes. In average-case or generic-case complexity, where the measure on the set of inputs matters, the equivalence between these two points of view needs to be verified. Most of our results are unchanged if we take $F(A)$ rather than $X^*$ as the set of inputs.

If $Y$ is another finite set of generators for $G$ and $\mathcal{D}_Y$ is the $Y$-version of the decision problem $\mathcal{D}$ then these two decision problem are equivalent from the point of view of worst-case complexity. Indeed, every generator $x \in X = A \cup A^{-1}$ can be written as a word in $F(Y)$. Thus every word in $X^*$ can be re-written in linear time as a word in $Y^*$ representing the same group element. This provides a linear-time reduction of $\mathcal{D}_A$ to $\mathcal{D}_Y$, and vice versa. Thus the worst-case complexity of group-theoretic decision problems does not depend on the choice of a finite generating set and is a true group invariant. By contrast, in the average or generic-case complexities a change in generating sets might conceivably give a different result and we will make invariance part of our definition. All of the results proved in this paper are invariant under change of a generating set.

A more complicated class of algorithmic problems can be described as *witness problems*. Unlike decision problems, a "witness problem" asks to produce, for a given element $u \in \mathcal{D}$, an explicit justification or "proof" of the fact that $u$ is, indeed, in $\mathcal{D}$. For example, the "witness" version of the Word Problem for a presentation $\langle A \mid R \rangle$, given a word $u \in ncl(R)$, asks for an explicit expression of $u$ as a product of conjugates of elements from $R^{\pm 1}$.

$$u = \prod_{i=1}^{t} u_j^{-1} r_j^{\epsilon_j} u_j,$$

where $u_j \in F(A)$, $r_j \in R$, and $\epsilon_j = \pm 1$.

The witness Conjugacy Problem would require producing a conjugating element, and the witness Membership Problem would ask for an expression of a given $u \in \langle v_1, \dots, v_k \rangle$ as a product of the given generators (and their inverses) of the subgroup considered. Although witness problems are increasingly important (for example, in group-based cryptography [3]), we concentrate here on the traditional decision problems.

Suppose we have a total algorithm $\Omega_1$ solving a decision problem $\mathcal{D}$ and also a partial algorithm $\Omega_2$ solving the problem generically with low generic-case complexity. Then by running $\Omega_1$ and $\Omega_2$ in parallel we obtain a new total algorithm $\Omega = \Omega_1 || \Omega_2$ which solves $\mathcal{D}$ with low generic-case complexity. The idea of putting these two algorithms together is in fact used by many practical experimenters. That is, for a particular problem one should look both for an exact solution with minimal known worst-case complexity and for a partial "generic" solution which will work very fast on most inputs. The computational group theory package "Magnus" already uses this philosophy very substantially, as most problems there are attacked by several algorithms running in parallel, including "fast checks" working with abelianizations and other quotients. We refer the reader to the article of G. Baumslag and C. F. Miller [13] for a more detailed discussion on "Magnus". More recently, several applications of genetic algorithms in group theory [49, 57] revealed that some classical problems that were believed to have only "too slow", i.e., non-practical, solutions, admit a very fast solution generically. This, as well as numerous computer experiments, provided an important source of intuition for the present paper.

If the generic-case complexity of $\Omega_2$ is very low and the worst-case complexity of the total algorithm $\Omega_1$ is not too high, then the combined algorithm may have low actual average-case complexity. The idea of using generic-case results to prove average-case results in this way seems very fruitful, and we have already been able to obtain some interesting results which will be the subject of a future paper.

## 3. Generic-case complexity

We have stressed that in order to measure the "largeness" of a set of words on an alphabet one needs a measure or, at least, an additive positive real-value function defined on some sets of words in the alphabet. For this paper we use the asymptotic density function suggested in the work of A. Borovik, A. Myasnikov and V. Shpilrain [17] (see also [16]) and similar in spirit to concepts considered by M. Gromov, A. Ol'shanskii and C. Champetier.

**Definition 3.1** (Asymptotic density). Let $X$ be a finite alphabet with at least two elements and let $(X^*)^k$ denote the set of all k-tuples of words on $X$. The *length* of a $k$-tuple $(w_1, ..., w_k)$ is the sum of the lengths of the $w_i$. Let $S$ be a subset of $(X^*)^k$. For every $n \geq 0$ let $B_n$ be the set of all $k$-tuples in $(X^*)^k$ of length at most $n$.

We define the *asymptotic density* $\rho(S)$ for $S$ in $(X^*)^k$ as

$$\rho(S) := \limsup_{n \to \infty} \rho_n(S)$$

where

$$\rho_n(S) := \frac{|S \cap B_n|}{|B_n|},$$

If the actual limit $\lim_{n \to \infty} \rho_n(S)$ exists, we denote $\widehat{\rho}(S) := \rho(S)$.

In the case where the limit

$$\lim_{n \to \infty} \rho_n(S) = \widehat{\rho}(S)$$

exists we shall be interested in estimating the speed of convergence of the sequence $\{\rho_n(S)\}$. To this end, if $a_n \geq 0$ and $\lim_{n \to \infty} a_n = 0$, we will say that the convergence is *exponentially fast* if there is $0 \leq \sigma < 1$ and $C > 0$ such that for every $n \geq 1$ we have $a_n \leq C\sigma^n$. Similarly, if $\lim_{n \to \infty} b_n = 1$ (where $0 \leq b_n \leq 1$), we will say that the convergence is *exponentially fast* if $1 - b_n$ converges to 0 exponentially fast.

**Definition 3.2** (Generic sets). We say that a subset $S \subseteq (X^*)^k$ is *generic* if $\widehat{\rho}(S) = 1$.

If in addition $\rho_n(S)$ converges to 1 exponentially fast, we say that $S$ is *strongly generic*.

What we have really defined is being *generic with respect to* $\widehat{\rho}$ in the sense discussed in Section 1. Since we now fix this particular concept of being generic, we simply say "generic" for the rest of this paper. The complement of a generic set is termed a *negligible* set. We can define *strongly negligible sets* in a similar manner. In the following lemma we collect several simple but useful properties of generic and negligible sets.

**Lemma 3.3.** *Let $S, T$ be subsets of $(X^*)^k$. Then the following hold:*
1) *$S$ is generic if and only if $\overline{S}$ is negligible.*
2) *If $S$ is generic and $S \subseteq T$ then $T$ is generic.*
3) *Finite unions and intersections of generic (negligible) sets are generic (negligible).*
4) *If $S$ is generic and $T$ is negligible, then $S - T$ is generic.*
5) *The collection $\mathcal{B}$ of all generic and all negligible sets forms an algebra of subsets of $(X^*)^k$.*

Now we can define generic-case complexity of algorithms.

**Definition 3.4** (Generic and strong generic performance of a partial algorithm). Let $\mathcal{D} \subseteq (X^*)^k$ be a decision problem and let $\mathcal{C}$ be a complexity class. Let $\Omega$ be a correct partial algorithm for $\mathcal{D}$, that is, whenever $\Omega$ reaches a definite decision on whether or not a tuple in $(X^*)^k$ belongs to $\mathcal{D}$, that decision is correct.

We say that $\Omega$ *solves $\mathcal{D}$ with generic-case complexity $\mathcal{C}$* if there is a generic subset $S \subseteq (X^*)^k$ such that for every tuple $\tau \in S$ the algorithm $\Omega$ terminates on the input $\tau$ within the complexity bound $\mathcal{C}$.

If in addition the set $S$ is strongly generic, then we say that the partial algorithm $\Omega$ *solves the problem $\mathcal{D}$ with generic-case complexity strongly $\mathcal{C}$*.

We again point out that we completely ignore the performance of $\Omega$ on tuples not in $S$ and the definition thus applies to the case where $\mathcal{D}$ has arbitrarily high worst-case complexity or is indeed undecidable.

One can now define "generic" complexity classes of decision problems in the obvious way.

**Definition 3.5** (Generic complexity classes)**.** Let $\mathcal{C}$ be a complexity class. Then Gen($\mathcal{C}$) denotes the class of all decision problems $\mathcal{D}$ for which there exists a partial algorithm solving $\mathcal{D}$ with generic-case complexity $\mathcal{C}$. Similarly, SGen($\mathcal{C}$) denotes the class of all decision problems $\mathcal{D}$ for which there exists a partial algorithm solving $\mathcal{D}$ with generic-case complexity strongly in $\mathcal{C}$.

As we mentioned before, while the worst-case complexity of most group-theoretic decision problems does not depend on the choice of a finite generating set for a group, it is not at all clear (and is probably false) that generic-case complexity *per se* is independent of the chosen set of generators. In order to have a true group-theoretic invariant, we need to incorporate such independence into the following definition.

**Definition 3.6** (Generic-case complexity of a decsiion problem **D** *for a group $G$*)**.** Let $G$ be a finitely generated group. Let **D** be an algorithmic problem about the group $G$. We say that *the decision problem* **D** *for $G$ has generic-case complexity in $\mathcal{C}$ (strongly in $\mathcal{C}$* if for *every* finite generating set $A$ of $G$ there exists a partial algorithm $\Omega(A)$ which solves the problem $\mathcal{D}_A \subset (A \cup A^{-1})^*$ with generic-case complexity $\mathcal{C}$ (strongly in $\mathcal{C}$).

## 4. Main results

In this section we formulate the main results of the paper. Even though our results regarding the word problem follow from the more general theorem about the membership problem (see Theorem B below), we state the word problem results first since most of the applications which we have in mind concern the word problem.

The concept of a group being non-amenable plays an important role in our results but for now the reader needs only to remember that any group which contains a free subgroup of rank two is non-amenable.

**Theorem A.** *Let $G$ be a finitely generated group. Suppose that $G$ has a finite index subgroup that possesses an infinite quotient group $\overline{G}$ for which the word problem is solvable in the complexity class $\mathcal{C}$. Then the word problem for $G$ has generic-case complexity in the class $\mathcal{C}$. Moreover, if the group $\overline{G}$* *is nonamenable, then the generic-case complexity of the word problem for $G$ is strongly in $\mathcal{C}$.*

There are a number of interesting immediate corollaries of the above result.

**Corollary 4.1.** *Let $G$ be a finitely generated group.*

1. *Suppose $G$ has a finite index subgroup that possesses an infinite word-hyperbolic quotient $\overline{G}$. Then the word problem for $G$ is generically in real time. Moreover, if $\overline{G}$ is non-elementary, then the word problem for $G$ is strongly generically in linear time.*
2. *Suppose $G$ has a finite index subgroup that possesses an infinite automatic quotient $\overline{G}$. Then the word problem for $G$ is generically in quadratic time. Moreover, if $\overline{G}$ is non-amenable, then the word problem for $G$ is strongly generically in quadratic time.*
3. *Suppose $G$ has a finite index subgroup that possesses an infinite quotient group $\overline{G}$, that is linear over a field of zero characteristic.*

   *Then the word problem for $G$ is generically in polynomial time. Moreover, if if $\overline{G}$ is not virtually solvable, then the word problem for $G$ is strongly generically in polynomial time.*

*Proof.* It is well known that for any word-hyperbolic group and for any finite generating set of this group, there is a set of defining relators for which Dehn's algorithm solves the word problem in linear time in the length of the input word. Moreover, this linear-time algorithm can be carried out by a multi-tape Turing machine. This was first observed by Domanski and Anshel [4] (see also [2] for a detailed description of the algorithm). Moreover, Holt and Rees [38, 39] have proved that for a word-hyperbolic group the algorithm solving the word problem can be carried out by a multi-tape *real-time* Turing machine.

It is also well known that any word-hyperbolic group is either virtually cyclic (in which case it is called *elementary*) or contains a free group of rank two (in which case it is called *non-elementary*). Thus every non-elementary word-hyperbolic group is non-amenable. Together with Theorem A this implies the first part of Corollary 4.1.

Similarly, the classical result of [27] shows that for an automatic group with any finite generating set there is an algorithm which solves the word problem in quadratic time. Again, by Theorem A the second part of Corollary 4.1 immediately follows.

An important result of R. Lipton and Y. Zalcstein [43] states that for a finitely generated group over a field of characteristic zero the word problem is solvable in log-space and hence in polynomial time. By a famous theorem of J.Tits [56], a finitely generated linear group that is not virtually solvable contains a nonabelian free subgroup and hence is nonamenable. This, together with Theorem A, implies the third part of Corollary 4.1. □

**Example 4.2.** If $G$ is any finitely generated group with infinite abelianization then $G$ maps onto the infinite cyclic group and hence by Corollary 4.1 the word problem in $G$ is solvable generically in linear time. This is also equivalent to being able to write $G$ as an HNN extension in some way. The result thus applies to all knot groups, all Artin groups and infinite one-relator groups.

**Example 4.3.** Let $G$ be a finitely generated infinite virtually solvable group. Then $G$ has a finite index subgroup that possesses an infinite virtually abelian quotient. Hence by Corollary 4.1 the word problem in $G$ is solvable generically in linear time.

**Example 4.4.** Recall that the $n$-strand braid group $B_n$, where $n \geq 3$, is given by the presentation

$$B_n = \langle a_1, \ldots, a_{n-1} \mid a_i a_{i+1} a_i = a_{i+1} a_i a_{i+1}, \text{ for } i = 1, ..., n-2,$$
$$\text{and } a_i a_j = a_j a_i \text{ for } |i-j| > 1 \rangle$$

The *pure braid group* $P_n$ corresponds to those $n$-strand braids where every strand ends in the same position that it begins. Then $P_n$ is a normal subgroup of index $n!$ in $B_n$ and $B_n/P_n$ is isomorphic to the symmetric group $S_n$. While it is hard to map one braid group onto another, this task is easy with pure braid groups: for $n \geq 4$ the group $P_n$ maps onto $P_{n-1}$ by "pulling out" the last strand of a braid. Thus for every $n \geq 3$ the group $P_3$ is a quotient group of $P_n$. It is well-known that $P_3 \cong F(a,b) \times \mathbb{Z}$. Thus for each $n \geq 3$ the group $P_n$ has a non-abelian free quotient $F(a,b)$. Since $P_n$ is of finite index in $B_n$ and since $F(a,b)$ is nonelementary word-hyperbolic, Corollary 4.1 implies that for $n \geq 3$ the groups $P_n$ and $B_n$ have word-problems solvable with generic-case complexity *strongly* in linear-time.

**Example 4.5.** Let $G = Aut(F_n)$ or $G = Out(F_n)$ where $n \geq 2$. Then by looking at the action of an automorphism (an outer automorphism) of $F_n$ on the abelianization of $F_n$ we see that $G$ maps onto the group $GL(n, \mathbb{Z})$. Since the word problem in $GL(n, \mathbb{Z})$ is solvable in quadratic time and $GL(n, \mathbb{Z})$ is non-amenable (it contains a nonabelian free subgroup), Corollary 4.1 implies hat the word problem for $G$ is solvable strongly generically in quadratic time.

This observation raises the interesting question of determining the worst-case complexity of the word problem for $G = Aut(F_n)$, say with generators the elementary Nielsen automorphisms. We usually think that the "obvious algorithm" of checking if the action of an automorphism $\alpha$ fixes the generators is very simple, but writing out all the intermediate steps could yield exponentially long words. It is not clear if there an algorithm with better worst-case complexity.

**Example 4.6.** Theorem A holds even if $G$ has unsolvable word problem. We consider the finitely presented Boone group $\mathcal{B}$ with unsolvable word problem as described in J. Rotman's book [52]. One proves the word problem unsolvable by showing that equality between certain "special" words exactly mimics the word problem in a semigroup with undecidable word problem. We again have the situation that the complexity hinges on words of a very special form. It is easy to see that the group $\mathcal{B}$ has the nonabelian free group generated by all the $r_i$ as the quotient group which is obtained by killing all the other generators. Thus the stronger conclusion of the theorem applies and the generic-case complexity of the word problem for $\mathcal{B}$ is strongly linear time. This is not really surprising and is a precise version of the statement that the group $\mathcal{B}$ is "large" and the set of special words is really quite "sparse".

**Example 4.7.** Let $G$ be a group with a finite presentation involving at least two more generators than relators. By the result of B. Baumslag and S. J. Pride [11] $G$ has a subgroup of finite index that can be mapped homomorphically onto the free group of rank two. Hence by Corollary 4.1 $G$ has word problem solvable strongly generically in linear time.

In strong contrast with worst-case complexity is the fact that generic-case complexity for a problem **D** for a group $G$ tells us nothing whatsoever about the complexity of **D** for subgroups of $G$. For example, if $G$ is any finitely generated group, then $G$ is certainly embedded in the direct product $P = G \times F(a,b)$ of $G$ and the free group $F(a,b)$ of rank two. We can apply Theorem A to $P$ by taking the homomorphism to $F(a,b)$ which kills all the elements of $G$. Since $F(a,b)$ is hyperbolic and non-amenable, Theorem A implies that the word problem in $P$ is strongly generically in linear time. But this says nothing at all about $G$ because we just erased all information about $G$. This remark does show that every finitely generated group can be embedded in a finitely generated group whose word problem has generic-case complexity strongly in linear time. A well-known theorem of B. H. Neumann (see [44]) shows that there are continuumly many 2-generator groups, and thus there are continuumly many $n$-generator groups for every $n \geq 2$. Thus there are continuumly many finitely generated groups whose word problem has generic-case complexity strongly linear time. This is in sharp contrast with the fact that there are only countably many finitely generated groups with solvable word problem.

The following computer experiment is easy to program. Let $F_n$ be a free group of rank $n$ and let $\phi$ be the homomorphism from $F_n$ to $F_{n-k}$ defined by sending the the first $k < n$ generators of $F_n$ to the identity. Pick a large length $l$ and use a random number generator to generate a large number of random freely reduced words of length $l$. If one calculates the ratio of the number of words $w$ with $\phi(w) \neq 1$ to the total of number of words generated, one observes exactly the phenomena predicted by the theory of random walks.

We now turn to the membership problem. It is necessary to discuss both a basic situation where the membership problem is solvable and also a basic result about undecidability of the membership problem. We first observe that if $G$ is any finitely generated group and $H$ is a subgroup of finite index, then the membership problem for $H$ in $G$ is decidable in linear time . Choose a finite set $A$ of generators of $G$. The *Schreier coset graph* $\Gamma(G, H, A)$ is defined as follows. The vertex set $V$ of $\Gamma(G, H, A)$ is the set of cosets $\{Hg | g \in G\}$. If $y \in A$ then there is an edge labeled by $y$ from $Hg$ to $Hgy$. Every edge in $\Gamma(G, H, A)$ with label $a \in A$ is equipped with a formal inverse edge labeled by $a^{-1}$. Thus $\Gamma(G, H, A)$ is an oriented labeled graph.

If $A$ is finite and $H$ has finite index in $G$ then the graph $\Gamma(G, H, A)$ is a finite. We can view $\Gamma(G, H, A)$ as the transition graph of a finite state automaton $M$ where the initial state and the only final state is the coset $H1 = H$. By the definition of the coset graph, for any word $w$ on the generators and their inverses, $M$ accepts $w$ if and only if $w \in H$. Thus the membership problem for $H$ is indeed decidable in linear time: given a word $w \in (A \cup A^{-1})^*$, read $w$ on the graph starting at the coset $H$ and see if one ends back at the coset $H$. A generalized version of these ideas is currently important in geometric group theory.

**Theorem B.** *Let $G$ be a finitely generated group and let $H \leq G$ be a finitely generated subgroup of infinite index. Let $G_1$ be a subgroup of finite index in $G$ such that $H \leq G_1$ and let $\phi : G_1 \to \overline{G}$ be an epimorphism. Assume that $\overline{H} = \phi(H)$ is contained in a subgroup $\overline{K}$ of infinite index in $\overline{G}$ and such that the membership problem for $\overline{K}$ in $\overline{G}$ is in the complexity class $\mathcal{C}$. Then the membership problem for $H$ in $G$ has generic-case complexity in $\mathcal{C}$. Moreover, if the Schreier coset graph $\Gamma(\overline{G}, \overline{K}, A)$ is non-amenable (for some and hence any finite generating set $A$ of $G$), then the generic-case complexity of the membership problem for $H$ in $G$ is strongly in $\mathcal{C}$.*

The "strong" conclusion of Theorem B holds, for example, if $\overline{G}$ is non-elementary hyperbolic group and $\overline{K}$ is a quasiconvex subgroup of $\overline{G}$. Indeed, in this case the coset graph $\Gamma(\overline{G}, \overline{K}, A)$ is non-amenable by a recent result of I. Kapovich [40]. Since the membership problem for a quasiconvex subgroup of a hyperbolic group is solvable in linear time, Theorem B implies that the membership problem for $H$ in $G$ is strongly generically in linear time.

**Example 4.8.** An Artin group is a group with a presentation

(1) $\qquad G = \langle a_1, \ldots, a_n | u_{ij} = u_{ji}, \text{ where } 1 \leq i < j \leq n, \rangle$

be an Artin group where for $i \neq j$

$$u_{ij} := \underbrace{a_i a_j a_i \ldots}_{m_{ij} \text{ times}}$$

The Coxeter group $C$ associated with $G$ is the quotient group obtained by setting the squares of the generators equal to the identity. In general, the membership problem may be unsolvable for a Coxeter group or an Artin group. A Coxeter group or an Artin group is of *extra-large type* if all $m_{ij} \geq 4$.

Any Coxeter group of extra-large type with at least three generators is a non-elementary hyperbolic group. Appel and Schupp[5] solved the membership problem for subgroups generated by subsets of the given generators in Artin groups of extra-large type, but very little is known about the membership problem for arbitrary finitely generated subgroups. If $H$ is a finitely generated subgroup of an Artin group $G$ such that the image $\overline{H}$ in the Coxeter quotient $C$ and is quasiconvex then the membership problem for $H$ in $G$ has generic-case complexity strongly linear time. Schupp [53] showed that all groups in a very extensive class of Coxeter groups are locally quasiconvex, that is, every finitely generated subgroup is quasiconvex. Also, one can check whether or not a finitely generated subgroup has infinite index in quadratic time. This provides a large set of examples of finitely generated subgroups of Artin groups where the generic-case complexity of the membership problem is strongly in linear time.

**Example 4.9.** A basic negative result about the membership problem is the theorem of K. Mihailova [46] that if $P_n = F_n \times F_n$ is the direct product of two copies of the free group $F_n$ of rank $n \geq 2$, then there are subgroups $H$ of $P_n$ with unsolvable membership problem (see [44]). Let

$$G = \langle x_1, ..., x_n | r_1, ..., r_m \rangle$$

be a finitely presented group with unsolvable word problem. By using the well-known Higman-Neumann-Neumann embedding of a finitely presented group into a 2-generator group, we may assume that $n$ is any integer which is at least 2. We use the ordered pair notation for elements of the direct product $P_n = F_n \times F_n$. Let $H$ be the subgroup of $P_n$ with generators

$$(x_1, x_1), ..., (x_n, x_n), (1, r_1), ..., (1, r_m)$$

Since the $r_i$ are defining relators for $G$, an easy argument shows that

$$(u, v) \in H \text{ if and only if } u = v \text{ in } G$$

Thus deciding membership in $H$ is equivalent to solving the word problem in $G$.

We point out that "genericity" is operating at three different levels when considering the membership problem. Let us fix $P_n$ as the direct product of two free groups of rank $n$. Call a subgroup $H$ a *subgroup of Mihailova type* if $H$ has a set of generators of the form (∗) above, which is very special. If we choose a random set of generators for a subgroup, it is very unlikely that they will be even close to being of Mihailova type. The remarks above showed that membership in a Mihailova subgroup $H$ is equivalent to the word problem for the group $G$ whose defining relators are the $r_i$. So just among subgroups of Mihailova type, if we choose the $r_i$ at random we encounter the phenomenon that finitely presented groups on a fixed set of generators are generically hyperbolic and thus the membership problem for the corresponding $H$ is still actually solvable in linear time. But Theorem B still applies to a particular Mihailova subgroup chosen to have unsolvable membership problem. All the known explicitly constructed groups with unsolvable word problem have at least infinite cyclic quotients, even after embedding into a two-generator group. That is, there is a homomorphism $\phi$ from $F_n$ to $\mathbb{Z}$ which sends all the $r_i$ to the identity. Let $\psi$ be the homomorphism from $P_n$ to $Q_n = F_n \times \mathbb{Z}$ defined by $\psi(u, v) = (u, \phi(v))$. The image $\overline{H}$ of $H$ is $F_n \times \{1\}$ which has infinite index in $Q_n$. The membership problem for $\overline{H}$ in $Q_n$ is clearly in linear time since to decide if $(u, v) \in \overline{H}$ one only has to check if $v$ equals the identity. If, for example, we use the Boone group $\mathcal{B}$ directly, without reducing the number of generators, to construct a Mihailova subgroup, then we have a homomorphism where the image $\overline{H}$ is the first factor of $F_k \times F_2$ and the generic-case complexity of the membership problem for $H$ is strongly linear time.

There is a similar theorem for the conjugacy problem.

**Theorem C.** *Let $G$ be a non-cyclic finitely generated group with infinite abelianization. Then the generic-case complexity of the conjugacy problem for $G$ is linear time.*

Theorem C is applicable to a wide variety of groups, such as infinite one-relator groups, Braid and Artin groups, knot groups etc.

We shall see that the proof of the theorem reduces to the case of the word problem since two words are conjugate in an abelian group if and only if they are equal. The reader has probably noticed that a statement about strong generic-case complexity in the case of non-amenable quotients is missing from the theorem. At the present writing we do not have a proof which is invariant under changing the set of generators although we believe that such a theorem is true.

A very interesting class of finitely presented groups with unsolvable conjugacy problem is the class of residually finite groups with unsolvable conjugacy problem constructed by C. F. Miller [47]. Given any finitely presented group $G$ with unsolvable word problem, Miller shows how to construct a group $M(G)$ which is the semidirect product of two finitely generated free groups (and which is thus residually finite) where conjugacy in $M(G)$ codes the word problem for $G$. As usual, the "code words" have a special form. The groups $M(G)$ have large nonabelian free quotients. We can show (although the argument is not presented in this article) that the conjugacy problem of such an $M(G)$ has generic-case complexity which is strongly linear time because the free quotient is obtained by simply killing some of the given generators.

We again stress some important limitations of generic case complexity.

First, just the definition of generic-case complexity does not say anything about the speed with which a particular sequence tends to one or zero. If the quotient group $\overline{G}$ is infinite but not "large enough", say $\overline{G} = \mathbb{Z}$, this speed may in fact be much slower than the exponentially fast convergence which we are really aiming at. The weaker convergence is all that we have for general one-relator groups.

Second, there is a substantial difference between our notion of "generic performance" and the notion of "average case complexity". In a situation like the word problem for one-relator groups where, although its complexity is not known, we at least have a total algorithm which is well understood, a future hope is to combine generic and worst-case methods to obtain average-case results. In this regard the work of [16, 17] about constructing explicit measures on free groups may be particularly useful.

In general, our approach simply shows that for the "decision" version of the word and the membership problem the fast " No" answer component of the set of all inputs is very large. One may be mainly interested in some infinite recursive subset of inputs and many examples may not admit algorithms with fast generic performance when restricted to the subset of interest.

Finally, our results do not say anything about the generic behavior of the "witness" versions of the word, conjugacy and membership problems. Yet it appears to us that if one has in mind practical cryptographic applications, these applications have to be based on the "witness" versions of the problems (rather than "decision" ones).

Thus we regard this paper as just the first step in the direction of understanding the generic-case and average-case behavior of various group-theoretic algorithms.

The results which we discuss in the last section of the paper (on finding languages which are *not* in given generic complexity classes) are due to Carl Jockusch and Frank Stephan. For example, the set of languages over a finite alphabet $A$ (with at least two letters) which are generically computable has measure zero (in a precise sense) in the set of all languages over $A$. Moreover, given any proper time-complexity function $f(n)$ one can construct a language that is deterministically computable in time $f^3(n)$ but which cannot be generically computed in time $f(n)$.

These general results do not, however, answer the question of constructing finitely generated *groups* with decision problems of arbitrarily "high" generic-case complexity, say with a word problem which is not generically solvable. All our results in this paper are proved by the "quotient method" of finding an infinite quotient group in which the relevant problems have the desired complexity. Using the existence of two disjoint recursively enumerable sets which are not recursively separable and the Adian-Rabin construction, C. F. Miller III [48] constructed an example of a finitely presented group $G$ all of whose nontrivial quotients have unsolvable word problem! This particular group $G$ therefore completely defeats our method of proof but it may well be the case that the word problem has low generic-case complexity for some different reason. Indeed, it seems to be a very difficult problem to construct a finitely generated group where the generic-case complexity of the word problem is provably not linear.

## 5. Cogrowth and simple random walks on regular graphs

The proofs of our theorems depend on already known nontrivial facts about the behavior of simple random walks on regular graphs. The really hard work is done by that theory, so we now turn to the results which we need.

The subject of random walks on graphs and groups is vast and very active. We refer the reader to [18, 33, 58, 63, 64] for some background information and further references in this area. We will recall several basic definitions in facts which are directly needed in our arguments.

**Definition 5.1.** Let $\Gamma$ be a $d$-regular graph (where $d \geq 2$) with a base-vertex $x_0$.

Then let $a_n(\Gamma) = a_n$ denote the number of reduced paths (i.e. paths without backtracks) of length $n$ from $x_0$ to $x_0$ in $\Gamma$. Similarly, let $b_n(\Gamma) = b_n$ be the number of all paths of length $n$ from $x_0$ to $x_0$ in $\Gamma$. Also let $r_n = r_n(\Gamma)$ be denote the number of reduced paths of length at most $n$ from $x_0$ to $x_0$ in $\Gamma$. Finally, let $z_n = z_n(\Gamma)$ denote the number of all paths of length at most $n$ from $x_0$ to $x_0$ in $\Gamma$. Thus $r_n = \sum_{i=0}^{n} a_i$ and $z_n = \sum_{i=0}^{n} b_i$.

Put

$$\alpha(\Gamma) = \alpha := \limsup_{n \to \infty} \sqrt[n]{a_n},$$

$$\beta(\Gamma) = \beta := \limsup_{n \to \infty} \sqrt[n]{b_n}$$

and

$$\nu(\Gamma) = \nu := \frac{1}{d}\beta(\Gamma).$$

We shall refer to $\alpha(\Gamma)$ as the *cogrowth rate* of $\Gamma$ and to $\nu(\Gamma)$ as the *spectral radius* of $\Gamma$. The number $\beta(\Gamma)$ will be called the *non-reduced cogrowth rate* of $\Gamma$.

It turns out that the definitions of $\alpha(\Gamma)$, $\beta(\Gamma)$ and $\nu(\Gamma)$ do not depend on the choice of a base-point $x_0 \in \Gamma$ and we have (see for example [64, 18]):

**Lemma 5.2.** *Let $\Gamma$ be a connected $d$-regular graph with a base-vertex $x_0$, where $d \geq 2$. Then:*
  1. *The values of $\alpha(\Gamma)$, $\beta(\Gamma)$ and $\nu(\Gamma)$ do not depend on the choice of a base-point $x_0 \in \Gamma$.*
  2. $0 \leq \alpha(\Gamma) \leq d-1$, $0 \leq \beta(\Gamma) \leq d$ and $0 \leq \nu(\Gamma) \leq 1$.
  3. $\nu = \limsup_{n \to \infty} \sqrt[n]{p^{(n)}}$ *where $p^{(n)}$ is the probability that a simple random walk on $\Gamma$ originating at $x_0$ will return to $x_0$ in $n$ steps.*

**Definition 5.3.** Let $\Gamma$ be a $d$-regular graph where $d \geq 2$. We will say that $\Gamma$ is *amenable* if $\nu(\Gamma) = 1$.

An important result connecting cogrowth and spectral radius was first obtained by R.Grigorchuck [33] and J.Cohen [24] for Cayley graphs of finitely generated groups and later generalized by L.Bartholdi [10] to the case of arbitrary regular graphs.

**Theorem 5.4.** Let $\Gamma$ be a $d$-regular graph (where $d \geq 2$). Let $\alpha = \alpha(\Gamma)$, $\beta = \beta(\Gamma)$ and $\nu = \nu(\Gamma)$.
Then

$$\nu = \begin{cases} \frac{\sqrt{d-1}}{d}\left(\frac{\alpha}{\sqrt{d-1}} + \frac{\sqrt{d-1}}{\alpha}\right) & \text{if } \alpha > \sqrt{d-1} \\ \frac{2\sqrt{d-1}}{d} & \text{otherwise.} \end{cases}$$

In particular $\nu < 1 \iff \alpha < d-1 \iff \beta < d$, that is $\nu = 1 \iff \alpha = d-1 \iff \beta = d$.

The above theorem indicates that $\Gamma$ is amenable if and only if it has maximal possible cogrowth for a $d$-regular graph.

The following classical result is known as Stolz' Theorem (see for example [55]):

**Lemma 5.5.** Suppose $x_n, y_n$ are sequences of real numbers such that $y_n < y_{n+1}$ for every $n$ with $\lim_{n \to \infty} y_n = \infty$ and such that a finite limit

$$\lim_{n \to \infty} \frac{x_{n+1} - x_n}{y_{n+1} - y_n}$$

exists. Then

$$\lim_{n \to \infty} \frac{x_{n+1} - x_n}{y_{n+1} - y_n} = \lim_{n \to \infty} \frac{x_n}{y_n}.$$

**Lemma 5.6.** Let $c_n \geq 0$ and $c > 1$ be such that $\lim_{n \to \infty} \frac{c_n}{c^n} = 0$. Put $f_n = \sum_{i=0}^n c_n$. Then $\lim_{n \to \infty} \frac{f_n}{c^n} = 0$

*Proof.* Applying Stolz' Theorem to $x_n = f_n$ and $y_n = c^n$ immediately yields Lemma 5.6. $\square$

Our principal technical tool is:

**Theorem 5.7.** Let $\Gamma$ be an infinite connected $d$-regular graph, where $d \geq 3$. Let $a_n = a_n(\Gamma)$ and $r_n = r_n(\Gamma)$. Then
(i)
$$\lim_{n \to \infty} \frac{a_n}{(d-1)^n} = \lim_{n \to \infty} \frac{b_n}{d^n} = 0.$$

(ii)
$$\lim_{n \to \infty} \frac{r_n}{(d-1)^n} = \lim_{n \to \infty} \frac{z_n}{d^n} = 0.$$

*Proof.* This fact is essentially due to L.Bartholdi as it follows from the remark on p.99 in [10]. It was first obtained (in a stronger form) by W.Woess [63] for the case where $\Gamma$ is the Cayley graph of a finitely generated group. We present briefly a formal argument for completeness.

Notice that by Lemma 5.6 (i) implies (ii) since $r_n = \sum_{i=0}^n a_n$ and $z_n = \sum_{i=0}^n b_n$. We will now verify (i).

Suppose first that $\alpha(\Gamma) < d - 1$ and hence $\beta(\Gamma) < d$ by Theorem 5.4. Then there is $N_0 \geq 1$ and $0 < a < d - 1$ such that for all $n \geq N_0$ we have $a_n \leq a^n$. Hence for $n \geq N_0$

$$\frac{a_n}{(d-1)^n} \leq \frac{a^n}{(d-1)^n} \xrightarrow[n \to \infty]{} 0$$

as required. A similar argument implies that $\lim_{n \to \infty} b_n/d^n = 0$. Hence the statement of Theorem 5.7 obviously holds. Thus we may assume that $\alpha(\Gamma) = d - 1$, so that $\beta(\Gamma) = d$ and $\nu(\Gamma) = 1$ by Theorem 5.4. Then the word-by-word repetition of the proof of Lemma 3.9 of [10] implies that

$$\lim_{n \to \infty} \frac{a_n}{(d-1)^n} = \lim_{n \to \infty} \frac{b_n}{d^n} = 0.$$

Indeed, Lemma 3.9 of [10] proves a stronger version of Theorem 5.7 under the assumption that $\Gamma$ is also quasi-transitive. However, the only place in the proof of Lemma 3.9 in [10], where quasi-transitivity is used, is to conclude that $\beta(\Gamma) = d$ which is already known in our case. □

In case where $\Gamma$ is non-amenable, we can say even more.

**Proposition 5.8.** *Let $\Gamma$ be a non-amenable connected $d$-regular graph where $d \geq 3$ (and hence $\Gamma$ is infinite). Let $a_n = a_n(\Gamma)$, $r_n = r_n(\Gamma)$, $b_n = b_n(\Gamma)$ and $z_n = z_n(\Gamma)$. Then*

1. *Both $\frac{a_n}{(d-1)^n} \to 0$ and $\frac{r_n}{(d-1)^n} \to 0$ exponentially fast.*
2. *Both $\frac{b_n}{d^n} \to 0$ and $\frac{z_n}{d^n} \to 0$ exponentially fast.*

*Proof.* Since $\Gamma$ is non-amenable, we have $\alpha = \limsup \sqrt[n]{a_n} < d - 1$ which immediately implies that $\frac{a_n}{(d-1)^n} \to 0$ exponentially fast. It also means that there are $n_0 \geq 1$ and $1 < a < d - 1$ such that for any $n \geq n_0$ we have $a_n \leq a^n$. Hence for $n \geq n_0$

$$r_n = r_{n_0-1} + \sum_{i=n_0}^n a_i \leq r_{n_0-1} + a^{n_0} \frac{a^{n-n_0} - 1}{a - 1}$$

Thus there are $A, B > 0$ such that for any $n \geq n_0$ we have $r_n \leq A + Ba^n$. Since $1 < a < d - 1$, this implies that $\frac{r_n}{(d-1)^n}$ converges to zero exponentially fast. Thus part 1 of Proposition 5.8 is verified.

The non-amenability of $\Gamma$ implies $\beta = \limsup \sqrt[n]{b_n} < d$, which implies part 2 of Proposition 5.8 by the same argument as above. □

## 6. Cogrowth in groups

Let $G$ be a group with a fixed finite generating set $A$ consisting of $k \geq 1$ elements. If $w$ is a word in $A \cup A^{-1}$, we will denote by $\pi(w)$ the element of $G$ represented by $w$. We will also denote by $|w|$ the length of the word $w$. For an element $g \in G$ denote by $|g|_A$ the length of a shortest word in $A \cup A^{-1}$ representing $g$. Also, if $Q$ is an alphabet, we will denote by $Q^*$ the set of all words in $Q$. For a subset $S \subseteq G$ we will denote by $S_A$ the set of all words in $(A \cup A^{-1})^*$ representing elements of $S$.

Let $H \leq G$ be a fixed subgroup (not necessarily normal). Let $\Gamma = \Gamma(G, H, A)$ be the *Schreier coset graph* defined in Section se:1-4. Then $\Gamma$ is a connected $2k$-regular graph. Note also that if $H$ is normal in $G$, then $\Gamma$ is precisely the Cayley graph of the group $G/H$ with respect to the generating set $A$. Thus every edge-path in $\Gamma$ has a label which is a word in the alphabet $A \cup A^{-1}$. It is easy to see that for any word $w$ and any vertex $x$ of $\Gamma$ there exists a unique path in $\Gamma$ with label $w$ and origin $x$. Moreover, if $w$ is the label of a path $p$ starting at the vertex $x_0 := H1$ in $\Gamma$, then $\overline{w} \in H$ if and only if the terminal vertex of $p$ is also equal to $H1$. The graph-theoretic concepts from the previous section can now be re-stated as follows:

$$a_n(G, H, A) = \#\{w | w \text{ is a freely reduced word of length } n \text{ in } A \cup A^{-1} \text{ with } \pi(w) \in H\},$$

$$b_n(G, H, A) = \#\{w | w \text{ is a word of length } n \text{ in } A \cup A^{-1} \text{ with } \pi(w) \in H\},$$

$$r_n(G, H, A) = \#\{w | w \text{ is a freely reduced word of length } \leq n \text{ in } A \cup A^{-1} \text{ with } \pi(w) \in H\}$$

and

$$z_n(G, H, A) = \#\{w | w \text{ is a word of length } \leq n \text{ in } A \cup A^{-1} \text{ with } \pi(w) \in H\}.$$

**Proposition 6.1.** *Let $G$ be a group with a fixed finite generating set $S$ and let $\Gamma = \Gamma(G, H, A)$ be the coset graph with base-vertex $x_0 = H1$. Then:*

$$a_n(G, H, A) = a_n(\Gamma), \quad b_n(G, H, A) = b_n(\Gamma), \quad r_n(G, H, A) = r_n(\Gamma) \quad \text{and}$$
$$z_n(G, H, A) = z_n(\Gamma).$$

*Proof.* This fact follows directly from the definition of $\Gamma = \Gamma(G, H, A)$ and the fact that a word $w$ over $A \cup A^{-1}$ represents an element of $H$ if and only if the path in $\Gamma$ staring at $H1$ and labeled $w$ terminates at $H1$.  □

For this reason $\alpha(G, H, A) := \alpha(\Gamma)$ is called the *co-growth rate* of $H$ in $G$ with respect to $A$ and $\beta(G, H, A) := \beta(\Gamma)$ is called the *non-reduced co-growth rate* of $H$ in $G$ with respect to $A$. Similarly, $\nu(G, H, A) := \nu(\Gamma)$ is called the *spectral radius* of $H$ in $G$ with respect to $A$. As before, $\alpha(G, H, A) \leq 2k-1$, $\beta(G, H, A) \leq 2k$. Moreover $\alpha(G, H, A) = 2k-1$ if and only if $\beta(G, H, A) = 2k$ if and only if $\Gamma$ is amenable.

It is easy to see (and it is well-known) that amenability of $\Gamma(G, H, A)$ does not depend on the choice of a finite generating set $A$ for $G$:

**Proposition 6.2.** *Let $G$ be a finitely generated group and $H \leq G$ be a subgroup. Suppose $A, B$ are two finite generating sets for $G$. Put $\Gamma = \Gamma(G, H, A)$ and $\Gamma' = \Gamma(G, H, B)$. Then $\Gamma$ is amenable if and only if $\Gamma'$ is amenable.*

*Proof.* By Proposition 38 and Theorem 51 of [18] amenability is a quasi-isometry invariant for regular graphs of finite degree. Let us equip $\Gamma$ and $\Gamma'$ with simplicial metrics $d$ and $d'$ accordingly. Let $C := \max\{|a|_B \,|\, a \in A\}$ and $C' := \max\{|b|_A \,|\, b \in B\}$. Then for any two cosets $Hg_1, Hg_2$ we have $d'(Hg_1, Hg_2) \leq C'd(Hg_1, Hg_2)$ and $d(Hg_1, Hg_2) \leq Cd'(Hg_1, Hg_2)$. Thus the identity map $Id : (V\Gamma, d) \to (V\Gamma', d')$ is a quasi-isometry, which implies the statement of the proposition. $\square$

According to the traditional definition, a finitely generated group $G$ is called *amenable* if any action of $G$ on a compact space $Y$ by homeomorphisms admits a $G$-invariant probability measure on $Y$. It turns out that if $A$ is finite generating set of $G$ and $H \leq G$ is normal, then $\Gamma = \Gamma(G, H, A)$ is amenable if and only if the quotient group $G_1 = G/H$ is amenable. In particular $G$ itself is amenable if and only if its Cayley graph $\Gamma(G, A)$ is amenable.

Suppose now that $G = F = F(x_1, \ldots, x_k)$ is a free group of rank $k \geq 2$. It is easy to see that the number of vertices of the $n$-sphere in the Cayley graph of $F$ with respect to the free basis $A = \{x_1, \ldots, x_k\}$ is $2k(2k-1)^{n-1}$ for $n \geq 1$. Hence the number of elements of $F$ in the $n$-ball around the identity is $1 + \frac{k}{k-1}[(2k-1)^{n-1} - 1]$ for $n \geq 1$.

**Theorem 6.3.** *Let $F = F(x_1, \ldots, x_k)$ and let $H \leq F$ be a subgroup, where $k \geq 2$. Put $A = \{x_1, \ldots, x_k\}$. Let $a_n = a_n(F, H, A)$, $r_n = r_n(F, H, A)$ and $\alpha = \alpha(F, H, A)$. Similarly, let $b_n = b_n(F, H, A)$, $z_n = z_n(F, H, A)$ and $\beta = \beta(F, H, A)$.*

*Then $\Gamma$ is a $2k$-regular graph, and $\alpha \leq 2k-1$, $\beta \leq 2k$. Moreover:*

1. *If $[F : H] = \infty$ then*

$$\lim_{n\to\infty} \frac{a_n}{(2k-1)^n} = \lim_{n\to\infty} \frac{r_n}{(2k-1)^n} = 0.$$

*and*

$$\lim_{n\to\infty} \frac{b_n}{(2k)^n} = \lim_{n\to\infty} \frac{z_n}{(2k)^n} = 0.$$

*Thus $H_A$ has zero asymptotic density in $(A \cup A^{-1})^*$ since the number of vertices in the ball of radius $n \geq 1$ in the Cayley graph $\Gamma(F, A)$ is $1 + \frac{k}{k-1}[(2k-1)^{n-1} - 1]$.*

2. *If the coset graph for $F$ relative $H$ is non-amenable (and hence $[F : H] = \infty$) then all the limits in part 1 converge to zero exponentially fast.*
3. *If $[F : H] < \infty$ then*

$$\limsup_{n \to \infty} \frac{a_n}{(2k-1)^n} > 0, \quad \limsup_{n \to \infty} \frac{r_n}{(2k-1)^n} > 0,$$

$$\limsup_{n \to \infty} \frac{b_n}{(2k)^n} > 0, \quad \limsup_{n \to \infty} \frac{z_n}{(2k)^n} > 0.$$

*Proof.* Parts 1 and 2 of this statement follows immediately from Theorem 5.7 and Proposition 5.8. We will now establish part 3 Theorem 6.3. Note that $r_n \geq a_n$ and $z_n \geq b_n$. Thus it suffices to check that lim-sups involving $a_n$ and $b_n$ are positive. Since $[F : H] < \infty$, there is a normal subgroup of finite index $H_1 \leq F$ such that $H_1 \leq H \leq F$. Then $a_n(F, H, A) \geq a_n(F, H_1, A)$ and $b_n(F, H, A) \geq b_n(F, H_1, A)$. So it suffices to consider the case where $H$ is normal of finite index $p$ in $F$. In this case the coset graph $\Gamma = \Gamma(F, H, A)$ is finite and has $p$ vertices. Thus $\Gamma$ is amenable and $\alpha(\Gamma) = 2k - 1$, $\beta(\Gamma) = 2k$. Then by the results of W.Woess [63] and L.Bartholdi [10]

$$\limsup_{n \to \infty} \frac{a_n}{(2k-1)^n} = \limsup_{n \to \infty} \frac{b_n}{(2k)^n} = \begin{cases} \frac{1}{p} & \text{if } \Gamma \text{ has some odd-length circuits} \\ \frac{2}{p} & \text{if } \Gamma \text{ has only even-length circuits} \end{cases}$$

Thus Theorem 6.3 is proved. □

When $H$ is a normal subgroup of $F$, the first part of Theorem 6.3 is originally due to W.Woess [63]. One can obtain much more precise statements than Theorem 6.3, where the denominators are replaced by powers of the co-growth rate of $H$, but Theorem 6.3 is quite sufficient for our purposes.

## 7. The Membership problem

We refer the reader to [2, 25, 34, 27, 31] for the background information on hyperbolic and automatic groups and their rational subgroups. We will recall several relevant definitions and results.

**Definition 7.1.** Let $G$ be a group with a finite generating set $A$. Let $L$ be a language over $A \cup A^{-1}$ such that $\pi(L) = G$, where $\pi$ is the natural map from the free semigroup on $A \cup A^{-1}$ to the group $G$. Let $H \leq G$ be a subgroup.

1. The subgroup $H \leq G$ is said to be *L-rational* if the set
$$L_H := \{w \in L \,|\, \pi(w) \in H\}$$
is a regular language and $H = \pi(L_H)$.
2. The subgroup $H \leq G$ is said to be *L-quasiconvex* if there exists $K > 0$ such that for any initial segment $u$ of a word $w \in L_H$ there is a word $v$ of length at most $K$ such that $\pi(uv) \in H$.

An important observation of S. Gersten and H. Short [31] states that:

**Proposition 7.2.** *Let $G$ be a group with a finite generating set $X$ and let $L$ be a language over $X \cup X^{-1}$ such that $\pi(L) = G$. Let $H \le G$ be a subgroup. Then $H$ is $L$-rational if and only if $H$ is $L$-quasiconvex.*

As the example of cyclic subgroups of $G = \mathbb{Z} \times \mathbb{Z}$ illustrates, it is possible that a particular subgroup is rational with respect to one automatic structure on $G$ but not the other. However, rationality is invariant in a somewhat weaker sense:

**Proposition 7.3.** *Let $G$ be an automatic group with a finite generating set $A$ and an automatic language $L$ over $A \cup A^{-1}$. Let $H \le G$ be an $L$-rational subgroup. Then for any finite generating set $B$ of $G$ there is an automatic language $L'$ over $B \cup B^{-1}$ for $G$ such that $H$ is $L'$-rational.*

*Suppose further that $G$ is word-hyperbolic. Then for any finite generating set $B$ of $G$ and for any automatic language $L'$ over $B \cup B^{-1}$ for $G$ the subgroup $H$ is $L'$-rational.*

*Proof.* The statement regarding hyperbolic groups is well-known and reflects the fact that for word-hyperbolic groups all possible notions of quasiconvexity for subgroups coincide.

The statement about automatic groups follows from the results of [27], although it is not stated there directly. Indeed, Theorem 2.4.1 of [27] proves that given $G, A, L$ as in Proposition 7.3, for any finite generating set $B$ of $G$ there is an automatic language $L'$ for $G$ over $B \cup B^{-1}$. The proof actually shows that any regular sub-language of $L$ gets "translated" into a regular sub-language of $L'$ with the same image in $G$. In this process $L_H$ gets "translated" in $L'_H$ and hence $L'_H$ is regular, as required. □

Because of Proposition 7.3 it is natural to adopt:

**Definition 7.4** (Rational Subgroup). Let $G$ be an automatic group and let $H \le G$ be a subgroup.

We say that $H$ is *rational* in $G$ if there exists an automatic language $L$ for $G$ such that $H$ in $L$-rational.

If $G$ is word-hyperbolic then a rational subgroup is also often referred to as *quasiconvex*.

**Proposition 7.5.** *Let $G$ be an automatic group and let $H \le G$ be a rational subgroup. Then:*
1. *For any finite generating set $A$ of $G$ there is an algorithm which solves the membership problem for $H$ in $G$ in quadratic time.*
2. *Suppose that $G$ is word-hyperbolic. Then for any finite generating set $A$ of $G$ there is an algorithm which solves the membership problem for $H$ is $G$ in linear time.*

*Proof.* Both of these statements are very well-known (see [27, 31, 29]), but we will indicate how the algorithm works.

To see (1) suppose that $A$ is a finite generating set of $G$. Then there is an automatic language $L$ over $A \cup A^{-1}$ for $G$ such that $L_H$ is regular. Given an arbitrary word $w$ over $A \cup A^{-1}$ we first apply the quadratic-time algorithm of [27] to take $w$ to a normal form in $L$, that is to find $w' \in L$ such that $w$ and $w'$ represent the same element of $G$. Since an automatic language $L$ consists of quasigeodesics [27], we have $|w'| \le c|w|$, where $c$ is some constant independent of $w$. We then check whether or not $w' \in L_H$ (which can be done in linear time in terms of $|w'|$). The total expanded time is clearly quadratic in $|w|$.

For a hyperbolic group $G$ and a rational subgroup $H \le G$ the algorithm solving the membership problem in linear time is virtually identical. First, for any finite generating set $A$ of $G$ there is a finite presentation of $G$ as $G = \langle A|R \rangle$ such that all Dehn-reduced words for this presentation are quasigeodesics (To see this one has to choose $R$ large enough and use the fact that local geodesics in the $\Gamma(G, A)$ are global quasigeodesics, provided the "local" parameter is chosen to be sufficiently large [2, 25, 32]). It is obvious that the set $L$ of all Dehn-reduced words is regular. Moreover, $H \le G$ is rational implies that $H$ is a quasiconvex subset of $\Gamma(G, A)$. Hence $H$ is $L$-quasiconvex since in a hyperbolic metric space a quasigeodesic and a geodesic with common endpoints are Hausdorff-close (again, see [2, 25, 32]). Therefore $H$ is $L$-rational by Proposition 7.2 and so $L_H$ is a regular language. Unlike the general case of an automatic group, as we mentioned earlier there is a *linear time* algorithm which takes a word $w$ over $A$ to its Dehn-reduced form $w'$ (see [2, 4]) where $|w'| \le |w|$. The algorithm solving the membership problem for $H$ in $G$ now works exactly as in the general automatic case. □

**Theorem B.** *Let $G$ be a finitely generated group and let $H \le G_1 \le G$ where $G_1$ has finite index in $G$. Suppose there is an epimorphism $\phi : G_1 \to \overline{G}$ such that $\overline{H} = \phi(H)$ is contained in a subgroup $\overline{K} \le \overline{G}$ of infinite index in $\overline{G}$ such that the membership problem for $\overline{K}$ in $\overline{G}$ is in the complexity class $\mathcal{C}$. Then the membership problem for $H$ in $G$ has generic-case complexity in $\mathcal{C}$. Moreover, if the coset graph of $\Gamma(\overline{G}, \overline{K})$ is non-amenable (for some and hence any finite generating set $A$ of $G$), then the generic-case complexity of the membership problem for $H$ in $G$ is strongly in $\mathcal{C}$.*

*Proof.* Let $A = \{x_1, \ldots, x_k\}$ be a finite generating set for $G$ and let $B$ be some finite generating set of $G_1$. Let $\pi : F \to G$ be the canonical epimorphism corresponding to the presentation $G = \langle x_1, \ldots, x_k \,|\, u_1, \ldots, u_m, \ldots \rangle$, where $F = F(x_1, \ldots, x_k)$. Let $K_1 := \phi^{-1}(\overline{K}) \le G_1 \le G$ and let $K_2 := \pi^{-1}(K_1) \le F$. Note that $[G_1 : K_1] = [\overline{G} : \overline{K}] = \infty$ and hence $[F : K_2] = \infty$. Moreover, the Schreier coset graphs $\Gamma(F, K_2, A) = \Gamma(G, K_1, A)$ and $\Gamma(G_1, K_1, B) = \Gamma(\overline{G}, \overline{K}, B)$ are quasi-isometric since $G_1$ has finite index in $G$. Thus $\Gamma(F, K_2, A)$ is non-amenable if and only if $\Gamma(\overline{G}, \overline{K}, B)$ is non-amenable.

Moreover $H \le K_1$. Thus if $w \in (A \cup A^{-1})^* - (K_2)_A$ then $\pi(w) \in G - H$. Let $z_n = z_n(F, K_2, A)$ let $C_n = \frac{(2k)^{n+1}-1}{2k-1}$ be the number of words in $(A \cup A^{-1})^*$ of length at most $n$.

Since $[F : K_2] = \infty$, Theorem 6.3 implies that $(K_2)_A$ has zero asymptotic density in $(A \cup A^{-1})^*$, that is

$$\lim_{n \to \infty} \frac{z_n}{C_n} = 0 \text{ and } \lim_{n \to \infty} \frac{C_n - z_n}{C_n} = 1,$$

and in both cases the convergence is exponentially fast if $\Gamma(\overline{G}, \overline{K}, B)$ is non-amenable. Thus the set $(A \cup A^{-1})^* - (K_2)_A$ is generic (and even strongly generic if $\Gamma(\overline{G}, \overline{K}, B)$ is non-amenable).

Fix a finite right Schreier transversal $T$ for $G_1$ in $G$ so that $1 \in T$, $|T| = [G : G_1]$ and $G = \cup_{t \in T} G_1 t$. Also fix the finite Schreier coset graph $\Gamma(G, G_1, A)$. Recall that a Schreier rewriting process for $G_1$ in $G$ consists in rewriting a word $w \in (A \cup A^{-1})^*$ as a word $vt$ where $v \in (B \cup B^{-1})^*$ and $t \in T$, so that $vt$ and $w$ represent the same element of $G$. Thus $w$ represents an element of $G_1$ if and only if $t = 1$. We recall, briefly, how the Schreier rewriting process works. For every $t \in T$ and $x \in A \cup A^{-1}$ we fix a word $u(t, x) \in (B \cup B^{-1})^*$ and an element $s(t, x) \in T$ such that $tx = u(t, x)s(t, x)$ in $G$. Given a word $w = x_1 \ldots x_n \in (A \cup A^{-1})^*$, where each $x_i \in A \cup A^{-1}$, we rewrite it as follows. First $1 \cdot x_1 = u(1, x_1)s(1, x_1)$. If $x_1 \ldots x_i$ has already been rewritten as $u_i t_i$, where $u_i \in (B \cup B^{-1})^*$ and $t_i \in T$, then

$$x_1 \ldots x_i x_{i+1} = u_i t_i x_{i+1} = u_i u(t_i, x_{i+1}) s(t_i, x_{i+1}).$$

Thus we put $u_{i+1} = u_i u(t_i, x_{i+1})$, $t_{i+1} = s(t_i, x_{i+1})$ and continue to the next step. At the end of the process we rewrite $w$ as $ut$, where $u$ is a word in $B \cup B^{-1}$ and $t \in T$. This rewriting process requires at most a linear amount of space and time in termsof $|w|$. Moreover, we have $|v| \le C|w|$, where $C = \max\{|u(t, x)| : t \in T, x \in A \cup A^{-1}\}$.

We will now construct a correct partial algorithm $\Omega$ for the membership problem of $H$ in $G$ as follows. Let $w$ be a word in $(A \cup A^{-1})^*$. Denote by $g$ the element of $G$ represented by $w$. First we read the word $w$ in the finite Schreier graph $\Gamma(G, G_1, A)$ starting from the vertex $G_1 \cdot 1$ and simultaneously apply the Schreier rewriting process to $w$. If the terminal vertex of the resulting path is different from $G_1 \cdot 1$, then $\pi(w) \notin G_1$ and hence $\pi(w) \notin H$. We declare that $w \notin MP(G, H, A)$ and terminate $\Omega$ in this case.

If the resulting path ends at $G_1 \cdot 1$ then $\pi(w) \in G_1$ and we have rewritten $w$ as a word $v$ in $(B \cup B^{-1})^*$. Note that $|v| \le C|w|$ where $C > 0$ is some constant independent of $w$.

By assumption the membership problem for $\overline{K}$ in $\overline{G}$ is solvable with complexity $\mathcal{C}$. We apply this algorithm to the word $v$. If the element $\overline{g}$ of $\overline{G}$ represented by $v$ does not belong to $\overline{K}$, then $\overline{g} \notin \overline{H}$. Hence the element $g$ of $G$ represented by $w$ and $v$ does not belong to $H$. In this case we declare that $w \notin MP(G, H, A)$ and terminate $\Omega$.

If it turns out that $v$ represents an element of $\overline{K}$, we terminate $\Omega$ without an answer.

The algorithm $\Omega$ terminates with a correct answer for every $w \notin (K_2)_A$. Since the set $(A \cup A^{-1})^* - (K_2)_A$ is generic (and even strongly generic if $\Gamma(\overline{G}, \overline{K}, B)$ is non-amenable), the statement of Theorem B holds. □

Our theorem on the word problem is an immediate corollary.

**Theorem A.** *Let $G = \langle x_1, ..., x_k | R \rangle$ be a finitely generated non-cyclic group. Suppose that $G$ has a finite index subgroup that possesses an infinite quotient group $\overline{G}$ in which the word problem is solvable in the class $\mathcal{C}$. Then the word problem for $G$ has generic-case complexity in the class $\mathcal{C}$.*

*Moreover, if the group $\overline{G}$ is non-amenable, then the generic-case complexity of the word problem for $G$ is strongly in $\mathcal{C}$*

*Proof.* Let $G_1 \leq G$ be a subgroup of finite index and let $\phi : G_1 \to \overline{G}$ be an epimorphism as in the statement of Theorem A. Put $H = \{1\} \leq G$ and $\overline{K} = \{1\} \leq \overline{G}$. Thus $\phi(H) \leq \overline{K}$. Moreover, the membership problem for $H$ in $G$ is precisely the word problem for $G$. Similarly the membership problem for $\overline{K}$ in $\overline{G}$ is precisely the word problem for $\overline{G}$. Now the conclusion of Theorem A follows from Theorem B. □

**Remark 7.6.** Theorem 6.3 shows that the statements of both Theorem A and Theorem B remain true if we define asymptotic density and genericity in terms of subsets of $F(A)$ (rather than subsets of $(A \cup A^{-1})^*$) by counting the ratios of the number of freely reduced words from a subset over the number of all freely reduced words.

**Corollary 7.7.** *Let $G$ be a finitely generated group and let $H \leq G_1 \leq G$, where $[G : G_1] < \infty$. Let $\phi : G_1 \to \overline{G}$ be an epimorphism with $\overline{H} = \phi(H)$. Then:*

1. *Suppose $\overline{G}$ is word-hyperbolic and $\overline{H} \leq \overline{G}$ is contained in a quasiconvex subgroup $\overline{K}$ of infinite index in $\overline{G}$. Then the membership problem for $H$ in $G$ is strongly generically in linear time.*
2. *Suppose $\overline{G}$ is automatic and $\overline{H} \leq \overline{G}$ is contained in a rational subgroup $\overline{K}$ of infinite index in $\overline{G}$. Then the membership problem for $H$ in $G$ is generically in quadratic time. Moreover, if $\Gamma(\overline{G}, \overline{K})$ is non-amenable then the membership problem for $H$ in $G$ is strongly generically in quadratic time.*

*Proof.* This follows directly from Theorem B and Proposition 7.5. □

## 8. The Conjugacy problem

Let $F = F(x_1, \ldots, x_k)$ and let $A = \{x_1, \ldots, x_k\}$ be a fixed free basis of $F$, where $k \geq 2$.

**Convention 8.1.** As before, we will denote by $C_n$ the number of words of length at most $n$ in $(A \cup A^{-1})^*$. Thus $C_n = \frac{(2k)^{n+1}-1}{2k-1}$.

Let $Q_n$ be the number of pairs $(w_1, w_2)$ of words in $(A \cup A^{-1})^*$ with $|w_1| + |w_2| \le n$.

Note that if $|w_1| + |w_2| = i \le n$ then $|w_1 w_2| = i \le n$. For a fixed word $w$ of length $i$ there are $(i+1)$ ways of representing $w$ as $w = w_1 w_2$. Recall that $A \cup A^{-1}$ consists of $2k$ letters.

Hence:

$$Q_n = \sum_{i=0}^{n} (i+1)(2k)^i$$

**Proposition 8.2.** *Let $H \le F$ be a subgroup of infinite index. Let $S \subseteq (A \cup A^{-1})^* \times (A \cup A^{-1})^*$ be the set of all pairs $(w_1, w_2)$ with $|w_1| + |w_2| \le n$ such that $w_1 w_2^{-1}$ represents an element of $H$. Then $\widehat{\rho}_A(S) = 0$.*

*Proof.* Let $b_j = b_j(F, H, A)$ be the number of all words of length $j$ representing elements of $H$. Then by Theorem 6.3 $\lim_{n \to \infty} b_j/(2k)^j = 0$ since $H$ has infinite index in $F$.

Suppose $(w_1.w_2)$ is a pair of words such that $|w_1| + |w_2| = i \le n$ and that the words $w := w_1 w_2^{-1}$ represents an element of $H$. For a fixed word $w$ of length $i$ representing an element of $H$ there are $i+1$ ways of writing $w$ as $w = w_1 w_2^{-1}$. Hence

$$\sigma_n(S) = \sum_{i=0}^{n}(i+1)b_i.$$

Therefore

$$\lim_{n \to \infty} \frac{\sigma_n(S)}{Q_n} = \lim_{n \to \infty} \frac{\sum_{i=0}^{n}(i+1)b_i}{\sum_{i=0}^{n}(i+1)(2k)^i} = \text{ (by Stoltz' Theorem)}$$
$$= \lim_{n \to \infty} \frac{(n+1)b_n}{(n+1)(2k)^n} = \lim_{n \to \infty} \frac{b_n}{(2k)^n} = 0,$$

as required. □

**Theorem C.** *Let $G$ be a non-cyclic finitely generated group with infinite abelianization. Then the generic-case complexity of the conjugacy problem for $G$ is linear time.*

*Proof.* Let $\overline{G}$ be the abelianization of $G$ and let $\phi : G \to \overline{G}$ be the abelianization map. Let $F = F(x_1, \ldots, x_k)$, $A = \{x_1, \ldots, x_k\}$ and let $\pi : F \to G$ be the presentation epimorphism. Let $H \le G$ be $H := Ker(\phi \circ \pi)$. As before, let $H_A$ be the set of all words in $(A \cup A^{-1})^*$ representing elements of $H$.

Let
$$S := \{(w_1, w_2) \in (A \cup A^{-1})^* \times (A \cup A^{-1})^* \,|\, \phi(\pi(w_1)) = \phi(\pi(w_2))\} =$$
$$= \{(w_1, w_2) \in (A \cup A^{-1})^* \times (A \cup A^{-1})^* \,|\, w_1 w_2^{-1} \in H_A\}.$$

By Proposition 8.2 $\widehat{\rho}_A(S) = 0$. If $(w_1, w_2) \notin S$ then $\phi(\pi(w_1)) \neq \phi(\pi(w_2))$ and hence $\phi(\pi(w_1))$ is not conjugate to $\phi(\pi(w_2))$ in $\overline{G}$ (since $\overline{G}$ is abelian). Thus if $(w_1, w_2) \notin S$ then $\pi(w_1)$ is not conjugate to $\pi(w_2)$ in $G$.

Since $\overline{G}$ is finitely generated abelian, there is an algorithm $\Omega$ which solves the word problem for $\overline{G}$ in linear time. Hence for any pair $(w_1, w_2) \notin S$ with $|w_1| + |w_2| \leq n$ the algorithm $\Omega$ will terminate in linear time of $n$ and declare that $\phi(\pi(w_1)) \neq \phi(\pi(w_2))$ and hence $\pi(w_1)$ is not conjugate to $\pi(w_2)$ in $G$. □

## 9. Some general observations on generic-case complexity

As mentioned in Sectionintro, we are greatly indebted to Carl Jockusch and Frank Stephan for stimulating conversations about some general features of generic-case complexity and the results in this section are due to them. First, Carl Jockusch observed that if we put a reasonable measure on the set of all languages over an alphabet $A$ with at least two letters, then the set of generically computable languages has measure zero. Second, Frank Stephan observed that the standard Time Hierarchy Theorem of complexity theory can be modified to separate deterministic time classes from generic complexity classes. Thus, for example, there is a language $L$ in DTIME($n^3$) which is *not* in GenTIME($n$).

Fix an alphabet $A$ with at least two letters. A language $L$ over $A$ is *generically computable* if there is a partial algorithm $\Omega$ such that the set $S$ on which $\Omega$ correctly decides membership in $L$ has $\widehat{\rho}(S) = 1$. The *canonical* or *shortlex* ordering of the set $A^*$ of all words on $A$ orders words first by length and within length, by the lexicographical ordering induced from a linear ordering of $A$. So we have a listing $\{w_1, ..., w_n....\}$ of $A^*$ in which all words of a shorter length come before all words of a longer length. We can now identify a language $L \subseteq A^*$ with its characteristic function $\chi_L$ where

$$\chi_L(n) = \begin{cases} 1 \text{ if } w_n \in L, \\ 0 \text{ if } w_n \notin L. \end{cases}$$

Since such a characteristic function is an infinite sequence $(b_n)_{n \geq 1}$ of 0's and 1's, we can regard it as the binary expansion of a real number in the unit interval $[0, 1]$. A binary expansion is unique except for those which are either all 0's or all 1's from some point onwards. A binary representation which is all 0's from some point onwards corresponds to a finite subset of $A^*$. There are only countably many finite subsets and excluding them gives a one-to-one correspondence between the infinite subsets of $A^*$ and the half-open interval $(0, 1]$. The standard Lebesgue measure on $(0, 1]$ then gives a measure on the set of infinite subsets of $A^*$ and this is the measure which we use.

**Theorem 9.1.** *Let $A$ be a finite alphabet with at least two letters. Fix a linear ordering of $A$ and let $m$ be the measure on the set of infinite languages over $A$ induced by the shortlex ordering as described above.*

*Then the set of languages over $A$ which are generically computable has measure zero.*

*Proof.* It suffices to show that if $\Omega$ is any fixed partial algorithm whose output is either 0 or 1 then the set of languages which are generically decided by $\Omega$ has measure 0. Since there are only countably many algorithms, it then follows that the set of all generically decidable languages has measure 0. Let $\omega$ be the infinite sequence of 0's and 1's where $\omega(n) = 1$ if $\Omega$ calculates 1 for $w_n \in A^*$ and $\omega(n) = 0$ otherwise. The point is that $\omega$ is now a *fixed* sequence.

For an integer $K \geq 1$ denote by $g(K)$ the the number of subsets of a set with $K$ elements which contain at least $3K/4$ elements of that set. We need only the fact that the ratio of $g(K)$ over the number $2^K$ of all subsets of a set with $K$ elements goes to 0 as $K \to \infty$. This follows easily from applying Stirling's formula and computing the asymptotics of the binomial coefficient $\binom{K}{3K/4}$. This computation shows that

$$\frac{\binom{K}{3K/4}}{2^K} = o(\sigma^K) \text{ as } K \to \infty$$

for some number $0 < \sigma < 1$. Hence

$$\frac{g(K)}{2^K} := \frac{\binom{K}{3K/4} + \binom{K}{3K/4+1} + \cdots + \binom{K}{K}}{2^K} \leq \frac{K}{4} \frac{\binom{K}{3K/4}}{2^K} \xrightarrow[K \to \infty]{} 0$$

For every integer $j \geq 0$ the set $A^*$ has exactly $s(j) := \frac{k^{j+1}-1}{k-1}$ words of length $\leq j$, where $k = \#A$. Thus the first $s(j)$ digits in the binary sequence of a language $L$ determine exactly which words in $A^*$ of length at most $j$ belong to $L$.

Fix an arbitrary $\epsilon > 0$. Take an integer $j_1 > 0$ large enough so that $\frac{g(K)}{2^K} \leq \frac{\epsilon}{2}$ for any integer $K \geq s(j_1)$. Let $Q_1$ be the set of all infinite binary sequences which agree with the first $s(j_1)$ digits of $\omega$ in at least $3s(j_1)/4$ positions. Note that for a fixed binary string $\alpha$ of length $s(j_1)$ the measure of the set of all infinite binary sequences with initial segment $\alpha$ is $2^{-s(j_1)}$. Hence $m(Q_1) \leq g(s(j_1))2^{-s(j_1)} \leq \frac{\epsilon}{2}$.

Now take an integer $j_2 > j_1$ large enough so that $\frac{g(K)}{2^K} \leq \frac{\epsilon}{2^2}$ for any integer $K > s(j_2)$. Let $Q_2$ be the set of all infinite binary sequences which agree with the first $s(j_2)$ digits of $\omega$ in at least $3s(j_2)/4$ positions. Again we see that $m(Q_2) \leq \frac{\epsilon}{2^2}$. Continue in this way, choosing at step $n$ an integer $j_n > j_{n-1}$ large enough so that for any integer $K \geq s(j_n)$ we have $\frac{g(K)}{2^K} \leq \frac{\epsilon}{2^n}$.

Let $Q_n$ be the set of all infinite binary sequences agreeing with the first $s(j_n)$ digits of $\omega$ in at least $3s(j_n)/4$ positions. Then $m(Q_n) \leq \frac{\epsilon}{2^n}$.

Put $Q = \cup_{n=1}^{\infty} Q_n$. Then

$$m(Q) \leq \sum_{n=1}^{\infty} \frac{\epsilon}{2^n} = \epsilon.$$

Now suppose that $L$ is any language generically decided by $\Omega$. Be our choice of the enumeration of $A^*$ and by the definition of generic computability, there exists an integer constant $i \geq 0$ such that for any $j \geq i$ the binary sequence of $L$ agrees with the initial segment of $\omega$ of length $s(j)$ in at least $3s(j)/4$ positions. Choose $n$ such that $j_n \geq i$. Then $\chi_L \in Q_n \subseteq Q$ by construction of $Q_n$.

Thus we have shown that for any $\epsilon > 0$ the set of all languages generically computable by $\Omega$ can be covered by a set of measure at most $\epsilon$. As required, this implies that the set of languages generically computable by $\Omega$ has measure zero. $\square$

The following theorem is due to Frank Stephan. Recall that we are following the definitions and notations of [51] for computational complexity. A *proper complexity function* $f$ is a non-decreasing function for which there is a multi-tape Turing machine which, on an input $w$ computes the string $1^{f(|w|)}$ in $O(|w| + f(|w|))$ steps and uses $O(f(|w|))$ space besides its input. The reason for insisting on proper complexity functions is that they can be used as "clocks" when simulating Turing machines. One effectively assigns a word $\gamma(M)$ on a fixed alphabet $A$ which codes the Turing machine $M$. There is a universal Turing machine $U$ which, when given as input a word $\gamma(M)w$, simulates the machine $M$ on the input $w$. (We can assume that $w$ is a word in the alphabet $\{0,1\}$.) If $f$ is a proper complexity function we can define a time-bounded version of the Halting Problem by

$$H(f) = \{\gamma(M)w : M \text{ accepts } w \text{ in at most } f(|w|) \text{ steps}\}.$$

The following statement is Lemma 7.1 of [51] which shows that, given the code of a Turing machine $M$, we do not need more than time $f^3(|w|)$ to simulate $M$ for $f(|w|)$ steps on an input $w$.

**Lemma 9.2.** $H(f) \in \text{DTIME}(f^3(n))$.

Using the lemma we can prove

**Theorem 9.3.** *If $f(n) \geq n$ is a proper complexity function then there is a language $L \subseteq \{0,1\}^*$ which is computable in time $f^3(n))$ but not generically computable in time $f(n)$.*

*Proof.* The idea of the proof is that for each Turing machine $M$ we specify infinitely many lengths devoted to "defeating" the machine $M$. We can do this by using ordered pairs. Let $\mathbb{N}^+$ denote the set of positive integers. The standard one-to-one enumeration, often called the "pairing function",

$$p : \mathbb{N}^+ \times \mathbb{N}^+ \to \mathbb{N}^+$$

is given by a simple formula and its inverse function $p^{-1}$ is also easily computable, certainly in cubic time. We define the language $L$ as follows. If $w$ is a word on $\{0,1\}$, let $n = |w|$ and calculate $p^{-1}(n) = (r, s)$. If $r$ is not the code $\gamma(M)$ of a Turing machine then $w \notin L$

If $r = \gamma(M)$ for some Turing machine $M$, we simulate the action of $M$ on the input $w$ for $f(|w|)$ steps. By Lemma 9.2 this requires at most $O(f^3(|w|))$ steps. Put $w$ in $L$ if and only if $M$ did not accept $w$ in $f(|w|)$ steps.

By construction we have $L \in \text{DTIME}(f^3(n))$. On the other hand, if $L$ were in $\text{GenTIME}(f(n))$ then there would exist a Turing machine $M'$ and an integer $n$ such that for all $m \geq n$ the machine $M'$ correctly decides membership in $L'$ on at least three-quarters of all words of length less than or equal to $m$. Let $r = \gamma(M')$, let $s > n$ and let $t = p^{-1}(r, s)$. Note that $t > n$. By construction, $M'$ does not correctly decide membership in $L'$ for any words of length $t$ in time $f(t)$. But more than half of the words of length less than or equal to $t$ have length exactly $t$. Hence $M'$ fails to generically decide $L'$ in time $f(n)$, yielding a contradiction. $\square$

Department of Mathematics, University of Illinois at Urbana-Champaign, 1409 West Green Street, Urbana, IL 61801, USA
*E-mail address*: kapovich@math.uiuc.edu

Department of Mathematics, The City College of New York, New York, NY 10031
*E-mail address*: alexeim@att.net

Department of Mathematics, University of Illinois at Urbana-Champaign, 1409 West Green Street, Urbana, IL 61801, USA
*E-mail address*: schupp@math.uiuc.edu

Department of Mathematics, The City College of New York, New York, NY 10031
*E-mail address*: shpil@groups.sci.ccny.cuny.edu